\documentclass[12pt,reqno]{article}
\usepackage{amsmath, amsthm, graphicx,amsfonts, amssymb,color}
\usepackage{bm}
\usepackage{hyperref}
\RequirePackage{amsthm,amsmath,amsfonts,amssymb,mathtools}
\RequirePackage[numbers]{natbib}

\usepackage{bm}
\usepackage{hyperref}
\usepackage{mathrsfs}

 \def\d{\mathrm{d}}\def\e{\mathrm{e}}
\newtheorem{theorem}{Theorem}
\newtheorem{lemma}[theorem]{Lemma}
 \newtheorem{corollary}[theorem]{Corollary}
 
 \newtheorem{remark}[theorem]{Remark}

 \def\<{\langle}\def\>{\rangle}
  \def\prf{\noindent{\bf Proof.~~}}
 \def\deprf{\hfill$\Box$\medskip}

\usepackage{mathrsfs}
\newcommand{\scr}[1]{\mathscr #1}

 \def\d{\mathrm{d}}

\def\de{\end{equation}}

\def\edar{\end{eqnarray}}

\def\l{\left}\def\r{\right}
\def\lan{\langle}\def\ran{\rangle}

\def\[{\l[} \def\]{\r]}
\def\({\l(} \def\){\r)}

\def\supp{{\rm supp\,}}

\def\roman1{{\uppercase\expandafter{\romannumeral1}}}
\def\roman2{{\uppercase\expandafter{\romannumeral2}}}

 \def\beqlb{\begin{eqnarray}}\def\eeqlb{\end{eqnarray}}
 \def\beqnn{\begin{eqnarray*}}\def\eeqnn{\end{eqnarray*}}
\def\d{{\mbox{\rm d}}}
 \def\supp{{\mbox{\rm supp}}}

\title{\bf  {Exponential and strong ergodicity for one-dimensional symmetric stable jump diffusions}}

\author{
	{\bf Tao Wang\footnote{
				Email: wang\_tao@mail.bnu.edu.cn}}\\
	\footnotesize{School of Mathematical Sciences, Beijing Normal University, }\\
	\footnotesize{Laboratory of Mathematics and Complex Systems, Ministry of Education}\\
	\footnotesize{Beijing 100875, China}\\
}

\date{}

\begin{document}

%%------------------------------------------------------------
 \maketitle
\begin{abstract}
	We obtain explicit criteria for both exponential ergodicity and strong ergodicity  for one-dimensional time-changed symmetric stable processes with  $\alpha\in(1,2)$.  Explicit lower bounds for  ergodic convergence rates are given.
\end{abstract}


%%------------------------------------------------------------
\section{Introduction and main results}
Let $X\coloneqq(X_t)_{t\geqslant 0}$ be a  symmetric $\alpha$-stable process on $\mathbb{R}$ with infinitesimal  generator $\Delta^{\alpha / 2}\coloneqq-(-\Delta)^{\alpha / 2}$, $\alpha\in (0,2)$, where $-(-\Delta)^{\alpha / 2}$ is the fractional Laplacian operator. It is well known that $X$ is  pointwise recurrent (i.e. it hits single
points almost surely) if and only if $\alpha \in (1,2)$ (see \cite[Remark 43.12]{Sa99}),  but it is not ergodic since the invariant measure is Lebesgue measure which is infinite.

Consider the following stochastic differential equation:
\begin{equation}\label{SDE}
	\mathrm{d} Y_{t}=\sigma\left(Y_{t-}\right) \mathrm{d} X_{t},
\end{equation}
where $\sigma$ is a strictly positive continuous function on $\mathbb{R}$.
By \cite[Proposition 2.1]{DK20}, there is a unique weak solution $Y=(Y_t)_{t\geqslant 0}$ to the SDE \eqref{SDE}, and $Y$ can also be expressed as a time change process $
Y_{t}\coloneqq X_{\zeta_{t}},$ {where}  $$\zeta_{t}\coloneqq\inf \left\{s > 0: \int_{0}^{s} \sigma\left(X_{u}\right)^{-\alpha} \mathrm{~d} u>t\right\}. 
$$
%and $\xi=\int_0^{\infty}\sigma(X_s)^{-\alpha}\d s$ is the explosive time. 
By %\cite[Theorem 5.2.2]{CM12} and 
\cite[Section 1.2]{CW14}, the generator of  $Y$ is  $L=\sigma^\alpha\Delta^{\alpha / 2}$  which is symmetric with respect to its invariant measure $\mu(\mathrm{d} x)=\sigma(x)^{-\alpha} \mathrm{d} x$.

Note that a time change does not change the  recurrence (cf.   \cite[Theorem 5.2.5]{CM12}). %but the reversible measure $\pi$ is changed. 
When $\alpha \in (1,2)$, $Y$ is pointwise recurrent, so that   %Thus $Y$ is Lebesgue irreducible.%and  $X$ is the one-dimensional symmetric stable process. 
it is Lebesgue irreducible (see \cite[Page 42]{DZ96} for the definition). 
%strong Feller property and stochastic continuity,
Thus by \cite[Proposition 4.1.1 and Theorem 4.2.1]{DZ96},    
$Y$ is {\it ergodic} whenever $\mu(\mathbb{R})<\infty$.
%So the study for various ergodic properties of  $Y$ is meaningful. 

Throughout this paper, we study explicit criteria  for both  exponential ergodicity and strong ergodicity for this process $Y$.  Furthermore,  
we obtain explicit estimates   for  ergodic convergence rates.

Now we assume that $\mu(\mathbb{R})<\infty$. Let $\pi(\d x)\coloneqq\mu(\d x)/\mu(\mathbb{R})$. By \cite[Section 1]{CW14},  the associated  regular Dirichlet form $(\mathscr{E},\mathscr{F})$ of $Y$ on $L^{2}(\pi)$ is given by
\begin{equation}\label{Diri form}
	\mathscr{E}(f, g)=\frac{1}{2} \int_{\mathbb{R}} \int_{\mathbb{R}}(f(x)-f(y))(g(x)-g(y)) \frac{C_{\alpha}\mathrm{d} x \mathrm{d} y}{|x-y|^{1+\alpha}} , \ \ f,g\in \mathscr{F},
\end{equation}
where 
\begin{equation}\label{Diri domain}
	\mathscr{F}=\{u\in L^2(\pi): \ \mathscr{E}(u,u)<\infty \}, 
\end{equation}
and  $C_\alpha=\frac{\alpha2^{\alpha-1}\Gamma((\alpha+1)/2)}{\sqrt{\pi}\Gamma(1-{\alpha}/{2})}$.

Denote by $(P_t)_{t\geqslant 0}$ the semigroup of $Y$,   $\pi(f)\coloneqq\int_{\mathbb{R}} f\d \pi$ and $\|f\|_{L^2{(\pi)}}\coloneqq\big(\pi(f^2))^{1/2}.$
We say that the   process $Y$  is {\bf exponentially ergodic (or $L^2$-exponentially convergent)}, if there exists $\lambda_1>0$, such that for any $f\in L^2(\pi)$,
\begin{equation}\label{def L2exp}
	\|P_tf-\pi(f)\|_{L^2(\pi)}\leqslant \e^{-\lambda_{1}t}\|f-\pi(f)\|_{L^2(\pi)},
\end{equation}
see \cite[(1.2)]{cmf98} for the definition and \cite[Page 158--160,(h) and (i)]{cmf05} for the equivalence between exponential ergodicity and $L^2$-exponential convergence.
The optimal constant $\lambda_1$ in  \eqref{def L2exp} is equal to the spectral gap
\begin{equation*}
	\lambda_{1}=\mathrm{gap}(\scr{E})\coloneqq\inf \{\mathscr{E}(f, f): f \in \mathscr{F}, \pi(f^{2})=1, \pi(f)=0\},
\end{equation*}
see \cite{cmf98} for more details.

Our first result is the explicit  criterion for exponential ergodicity and the explicit estimate for $\lambda_{1}$.
\begin{theorem}[Exponential ergodicity]\label{exp erg}
	$Y$ is exponentially ergodic
	if and only if $$\delta\coloneqq\sup_{x} |x|^{\alpha-1}\int_{\mathbb{R}\setminus (-|x|,|x|)}\sigma(y)^{-\alpha}\d y<\infty.$$
	Furthermore, 
	$$\lambda_1	\geqslant
	\frac{1}{4\omega_{\alpha}\delta},$$
	where 
	\begin{equation}\label{omegaalpha}
		\omega_{\alpha}\coloneqq-\frac{1}{\cos(\pi\alpha/2)\Gamma(\alpha)}>0.
	\end{equation}
\end{theorem}

For any open set $B\subset\mathbb{R}$, define the local Dirichlet eigenvalue by 
\begin{equation}\label{vari diri}
	\lambda_{0}(B)=\inf \{\mathscr{E}(f, f): f \in \mathscr{F}, \pi(f^{2})=1 \text { and }f|_{B^{c}}=0\}.
\end{equation}
In particular, denote by $\lambda_0\coloneqq\lambda_{0}(\{0\}^c)$. 

$\lambda_{0}(B)$ is the bottom of spectrum for the part Dirichlet form $(\scr{E}, \scr{F}^B)$ (see Section 2 for more detail). The probabilistic meaning of $\lambda_{0}(B)$ is 
the $L^2$-decay rate for the killed semigroup $P_t^B$ (see Section \ref{auxiliary} for the definition), i.e. 
\begin{equation*}
	\|P_t^{B}g\|_{L^2(\pi)}%=\int_{\mathbb{R}}\left(\int_{0}^\infty\e^{-\lambda t}\d E_{\lambda}g\right)^2\d \pi
	\leqslant \e^{-\lambda_0(B) t}\|g\|_{L^2(\pi)},\quad \text{for any}\ g\in L^2(\pi).
\end{equation*}
$\lambda_{0}(B)$ and $\lambda_1$ are closely related. By \cite[Theorem 1.4]{cmf00} or \cite[Theorem 1.1]{wfy99}, we know that if $\pi(B^c)>0$, then 
$\lambda_{1}\leqslant{\lambda_{0}(B)}/{\pi(B^c)}$. According to \cite[Proposition 3.2]{cmf00'},  $\lambda_{1}\geqslant\lambda_{0}$. Indeed, the sufficiency of Theorem \ref{exp erg} is based on the following result.
%The following theorem gives the estimate of $\lambda_{0}(B)$ when $B=\mathbb{R}\setminus\{0\}$ or $(0,\infty)$.
\begin{theorem}[Dirichlet eigenvalues]\label{Diri eigen}
	\ \ \ 
	
	(1) If $$\delta_{+}\coloneqq\sup _{x>0} x^{\alpha-1} \int_{x}^{\infty} \sigma(z)^{-\alpha} \mathrm{d} z<\infty,$$
	then		$$
	\lambda_{0}((0, \infty)) \geqslant \frac{(\alpha-1)\Gamma(\alpha / 2)^{2}}{4\delta_+}.	$$
	
	(2) If $$\delta\coloneqq\sup_{x} |x|^{\alpha-1}\int_{\mathbb{R}\setminus (-|x|,|x|)}\sigma(y)^{-\alpha}\d y<\infty,$$
	then 
	$$
	\frac{2}{\omega_{\alpha}}\left(\frac{1}{\delta_+}+\frac{1}{\delta_-}\right)\geqslant\lambda_{0}\geqslant 
	\frac{1}{4\omega_{\alpha}\delta},$$
	where $\omega_\alpha$ is given by \eqref{omegaalpha} and
	$$\delta_{-}\coloneqq\sup _{x>0} x^{\alpha-1} \int_{-\infty}^{-x} \sigma(z)^{-\alpha} \mathrm{d} z.%\quad \omega_{\alpha}\coloneqq-\frac{1}{\cos(\pi\alpha/2)\Gamma(\alpha)}>0.
	$$
\end{theorem}

Next, we study the strong ergodicity for $Y$. For this, 
let $\|\nu\|_{\rm{Var}}\coloneqq\sup_{|f|\leqslant 1}|\nu(f)|$ be the total variation of a signed measure $\nu$, and $P_t(x,\cdot)$ be the transition function. 
We say that	$Y$ is {\bf strongly ergodic}, if there exist  constants $1<C<\infty$ and $\kappa>0$, such that
$$\sup_{x\in\mathbb{R}}\|P_t(x,\cdot)-\pi\|_{\rm{Var}}\leqslant C\e^{-\kappa t}.$$
The optimal convergence rate 
$$\kappa= -\lim\limits_{t\rightarrow \infty}\frac{1}{t}\log \sup\limits_{x\in \mathbb{R}}\|P_t(x,\cdot)-\pi\|_{\mathrm{Var}},$$
see \cite{myh06} for more details.
\begin{theorem}[Strong ergodicity]\label{ergodc}
	$Y$ is strongly ergodic
	if and only if
	\begin{equation}\label{I}
		I\coloneqq\int_{\mathbb{R}}\sigma(x)^{-\alpha}|x|^{\alpha-1}\d x<\infty.
	\end{equation}
	Moreover, the optional convergence rate  $\kappa$
	in the strong ergodicity satisfies%(see \cite{myh06} for more details)
	$$\kappa\geqslant\frac{1}{\omega_{\alpha}I}>0.$$ 
\end{theorem}

\begin{remark}
	%	(1) Note that a symmetric $\alpha$-stable process on $\mathbb{R}^{d}$ is recurrent if and only if $\alpha\geqslant d$, i.e. $d=1$ and $\alpha\in [1,2)$. Therefore, we only need to consider the ergodicity for one-dimensional case.
	
	(1) Note that when $\alpha=1$, the process $Y$ is neighborhood recurrent but not pointwise recurrent 
	(see \cite[Section 3.3]{DK20}; also see \cite[Theorem I.1.5]{KAE18} for the criteria of general L\'{e}vy processes). Therefore, in this case, the Green operator $U^{(0)}$ killed on hitting the origin, which is important for the proofs of Theorem \ref{exp erg} and Theorem \ref{ergodc}, is not valid, for example, $U^{(0)}\mathbf{1}=\infty$. 
	
	(2) Our criteria are somehow comparable with those in the special case  $\alpha= 2$, the time-changed Brownian motion  on $\mathbb{R}$. 
	In the latter, the process is exponentially ergodic if and only if $$\delta_1\coloneqq\sup _{x>0} x \int_{\mathbb{R} \backslash(-x, x)} \sigma(z)^{-2} \mathrm{d} z<\infty,$$ and 
	$\lambda_{1}\geqslant\lambda_{0}\geqslant(4\delta_1)^{-1}$, 
	while it is  strongly ergodic if and only if	$$\int_{\mathbb{R}}\sigma(x)^{-2}|x|^{\alpha-1}\d x<\infty.$$ See \cite[Table 5.1 and Theorem 5.8]{cmf05}.
	
	(3)	By 	\cite[Table 2]{DK20}, for $\alpha>1$,  $\pm \infty$ is entrance boundary if and only if  \eqref{I} holds.  Therefore, Theorem \ref{ergodc} indicates that for a pointwise recurrent  time-changed symmetric stable process, the strong ergodicity is equivalent to entrance at $\pm\infty$. 
\end{remark}
By using Lyapunov functions, \cite[Theorem 1.7]{CW14} obtained some sufficient conditions for exponential ergodicity and  strong ergodicity, which now can be  derived by our Theorems \ref{exp erg} and \ref{ergodc}.

\begin{corollary}\label{lyapunov}
	$Y$ is exponentially ergodic if
	\begin{equation}\label{str lya}
		A_1\coloneqq\liminf _{|x| \rightarrow \infty} \frac{\sigma(x)}{|x|}>0,
	\end{equation}
	and $Y$ is strongly ergodic if
	\begin{equation}\label{exp lya}
		A_2\coloneqq\liminf _{|x| \rightarrow \infty} \frac{\sigma(x)}{|x|^{\gamma}}>0
	\end{equation}
	for some constant $\gamma>1$.
\end{corollary}  
\prf
From \eqref{str lya},  there exists $N_1>0$ such that for any $x>N_1$,
$\sigma(x)\geqslant A_1|x|/2$, so that  $$\delta\approx\sup_{|x|>N_1}|x|^{\alpha-1}\int_{\mathbb{R}\setminus(-|x|,|x|)}\sigma(y)^{-\alpha}\d y\leqslant \frac{2^{\alpha}}{(\alpha-1)A_1^{\alpha}}<\infty,$$
where the symbol ``$A\approx B$'' means that there exists $0<c_1,c_2<\infty$, such that $c_1 B\leqslant A\leqslant c_2B$.
Then  $Y$ is exponentially ergodic by Theorem \ref{exp erg}.

From \eqref{exp lya},  there exists $N_2>0$ such that for any $|x|>N_2$, $\sigma(x)\geqslant A_2|x|^\gamma/2,$ so that  $$I\approx\int_{[-N_2,N_2]^c}\sigma(x)^{-\alpha}|x|^{\alpha-1}\d x\leqslant %\left(\frac{A_2}{2}\right)^{-\alpha}\int_{[-N_2,N_2]^c}|x|^{\alpha(1-\gamma)-1}\d x=
\frac{2^{\alpha}N_2^{\alpha(1-\gamma)}}{\alpha(\gamma-1)A_2^{\alpha}}<\infty.$$
Then $Y$ is strongly ergodic by Theorem \ref{ergodc}.
\deprf

Actually, if $\sigma$ is a polynomial, then we have explicit results for the ergodicity and the convergence rates of $Y$. %In particular, the following estimate for $\lambda_0$ is sharp. 
\begin{corollary}\label{polyno}
	Consider the polynomial case:  $\sigma(x)=(1+|x|)^{\gamma}$. 
	
	(1) $Y$ is ergodic if and only if $\gamma>1/\alpha$.
	
	(2) For $\gamma>1$,
	$$	\frac{2(\alpha\gamma-1)^{\alpha\gamma}(\alpha(\gamma-1))^{\alpha(1-\gamma)}}{\omega_{\alpha}(\alpha-1)^{\alpha-1}}\geqslant\lambda_{0}\geqslant 	\frac{(\alpha\gamma-1)^{\alpha\gamma}(\alpha(\gamma-1))^{\alpha(1-\gamma)}}{8\omega_{\alpha}(\alpha-1)^{\alpha-1}},$$ 
	and for $\gamma=1,$ $$\frac{2(\alpha-1)}{\omega_{\alpha}}\geqslant\lambda_{0}\geqslant\frac{\alpha-1}{8\omega_{\alpha}}.$$ 
	
	(3) $Y$ is exponentially ergodic if and only if $\gamma\geqslant1$. Moreover, $$\lambda_{1}\geqslant\frac{(\alpha\gamma-1)^{\alpha\gamma}(\alpha(\gamma-1))^{\alpha(1-\gamma)}}{8\omega_{\alpha}(\alpha-1)^{\alpha-1}}$$ for $\gamma>1$ and  $\lambda_{1}\geqslant(\alpha-1)/{8\omega_{\alpha}}$ for $\gamma=1$. 
	
	(4) $Y$ is strongly ergodic if and only if $\gamma>1$. Furthermore, $\kappa\geqslant\alpha(\gamma-1)/2\omega_\alpha$.

	%for $\gamma>1$, $\lambda_{1}=\kappa\geqslant\lambda_0$.
\end{corollary}
\prf	
(1) Note that $\mu(\mathbb{R})=\int_{\mathbb{R}}(1+|x|)^{-\alpha\gamma}\d x<\infty$ if and only if $\gamma>1/\alpha$.  

(2) For $\gamma>1$,  a direct calculation shows that $\delta=2\delta_+=2\delta_-$ and 
\begin{equation}\label{gamma>1}
	\delta_+=\sup_{x} |x|^{\alpha-1}\int_{|x|}^{\infty}(y+1)^{-\alpha\gamma}\d y=\frac{(\alpha-1)^{\alpha-1}}{(\alpha\gamma-1)^{\alpha\gamma}(\alpha(\gamma-1))^{\alpha(1-\gamma)}}>0.
\end{equation}
Hence by Theorem \ref{Diri eigen},	
$$	\frac{2(\alpha\gamma-1)^{\alpha\gamma}(\alpha(\gamma-1))^{\alpha(1-\gamma)}}{\omega_{\alpha}(\alpha-1)^{\alpha-1}}\geqslant\lambda_{0}\geqslant 	\frac{(\alpha\gamma-1)^{\alpha\gamma}(\alpha(\gamma-1))^{\alpha(1-\gamma)}}{8\omega_{\alpha}(\alpha-1)^{\alpha-1}}.$$ 

For $\gamma=1$, we have $\delta=2\delta_+=2\delta_-$ and 
\begin{equation}\label{gamma=1}
	\delta_+=\sup_{x} |x|^{\alpha-1}\int_{|x|}^{\infty}(y+1)^{-\alpha}\d y=\lim_{|x|\rightarrow\infty}\frac{(|x|+1)^{1-\alpha}|x|^{\alpha-1}}{\alpha-1}=\frac{1}{\alpha-1}>0.
\end{equation}
So by Theorem \ref{Diri eigen},	\begin{equation*}
	\frac{2(\alpha-1)}{\omega_{\alpha}}\geqslant\lambda_{0}\geqslant\frac{\alpha-1}{8\omega_{\alpha}}.
\end{equation*}
%and the spectral gap $$\lambda_1	\geqslant\frac{\alpha-1}{8\omega_{\alpha}}.$$

(3) By \eqref{gamma>1}, \eqref{gamma=1} and Theorem \ref{exp erg}, we obtain the lower bound for $\lambda_1$, and see that if $\gamma\geqslant1$, $Y$ is exponentially ergodic.  

If $1/\alpha<\gamma<1$, then 
\begin{equation*}
	\delta=\sup_{x} 2|x|^{\alpha-1}\int_{|x|}^{\infty}(y+1)^{-\alpha\gamma}\d y=\lim_{|x| \rightarrow \infty}\frac{2}{\alpha\gamma-1}|x|^{\alpha-1}(|x|+1)^{1-\alpha\gamma}=\infty,
\end{equation*}
hence by Theorem \ref{exp erg}, $Y$ is exponentially ergodic if and only if $\gamma\geqslant1$. 

(4)	Note that  for  $1/\alpha<\gamma\leqslant1$,
\begin{equation*}
	\begin{split}
		I&\geqslant2\int_0^1\frac{x^{\alpha-1}}{(1+x)^{\alpha\gamma}}\d x+2\int_1^\infty\frac{x^{\alpha-1}}{(1+x)^{\alpha\gamma}}\d x\\
		&\geqslant\frac{1}{\alpha2^{\alpha\gamma-1}}+\frac{1}{2^{\alpha\gamma-1}}\int_{1}^{\infty}x^{\alpha-1-\alpha\gamma}\d x=\infty,
	\end{split}
\end{equation*} while  for $\gamma>1$, $I\leqslant{2}{\alpha^{-1}(\gamma-1)}^{-1}<\infty$. 
Thus Theorem \ref{ergodc} gives that $Y$ is strongly ergodic if and only if $\gamma>1$, and $\kappa\geqslant(\omega_{\alpha}I)^{-1}\geqslant\alpha(\gamma-1)/2\omega_\alpha$.

\deprf

\section{Killed process, Green function and time change}\label{auxiliary}
We first  recall some   definitions and properties.
Given an open set $B \subset \mathbb{R}$, denote by $$\tau_{B}\coloneqq\inf \left\{t > 0: Y_{t} \notin B\right\}$$  the {\it first exit time} from $B$ of the time-changed symmetric stable process $Y$.
Let $Y^{B}$ be the sub-process of $Y$ killed upon leaving $B$, whose transition function is
$$
P_{t}^{B} (x,A)\coloneqq\mathbb{P}_{x}\left[Y_{t}\in A, t<\tau_{B}\right],\ \text{for any}\ x\in B, \ \text{and Borel set}\ A\subset \mathbb{R}.$$
The  Green potential measure of the killed process $Y^B$ starting from $x$ is a Borel measure defined by
$$U^B(x,\d y)\coloneqq\int_0^{\infty}P_t^B(x,\d y)\d t .$$
The Green operator $U^B$ is given by
$$U^Bf(x)\coloneqq\int_Bf(y)U^B(x,\d y),\ x\in\mathbb{R},$$
for $f\in \mathscr{B}(\mathbb{R})$ with $U^B|f|<\infty$.

Recall that $(\scr{E}, \scr{F})$ is the Dirichlet form of $Y$ given by \eqref{Diri form} and \eqref{Diri domain}. Denote by 
$(\scr{E}, \scr{F}^B)$ the part Dirichlet form, where %associated with the killed process $X^B$:
$$ \mathscr{F}^B\coloneqq\{f\in \mathscr{F}, \widetilde{f}=0,\ \text{q.e. on}\ B^c \},$$
q.e. stands for quasi-everywhere, %i.e. it holds outside a set having zero 1-capacity with respect to  $(\scr{E}, \scr{F})$; 
and $\widetilde{f}$ is a quasi-continuous modification of $f$ 
(cf. \cite[Section 2.2]{OY13}).

%If  $(\scr{E}, \scr{F})$ is regular, then $(\scr{E}, \scr{F}^B)$ is a regular Dirichlet form on $L^2(\pi(\cdot\cap B)/\mu(B))$ (see  \cite[Theorem 3.5.7]{OY13}).
By \cite[Theorem 3.5.7]{OY13}, $(\scr{E}, \scr{F}^B)$ is a symmetric regular Dirichlet form, and  $Y^B$ is the process associated with  $(\scr{E}, \scr{F}^B)$.  %$\mathbf{1}_B\pi(\d x)$. %(cf. see \cite[section 2]{LLS11}). 

Note that for any nonempty set $B$,   $(\scr{E}, \scr{F}^B)$ is a transient Dirichlet form.
By \cite[Theorem 1.3.9]{OY13}, for any $f$ with $\int_{\mathbb{R}}|f(x)|U^B|f|(x) \pi(\d x)<\infty$, we have $U^Bf\in \scr{F}^B$, and for any $u\in\scr{F}^B$,
\begin{equation}\label{poisson weak-solution}
	\scr{E}(U^Bf,u)=\int fu\d \pi.
\end{equation}

Let  $\lambda_{0}(B)$ be  the bottom of spectrum for $(\scr{E}, \scr{F}^B)$: 
$$
\lambda_{0}(B)\coloneqq\inf \{\mathscr{E}(f, f): f \in \scr{F}^B, \pi(f^{2})=1 \}.$$
Denote by $C_0(B)$ the space of  continuous functions with compact support on $B$. Since $\scr{F}^B$ is the closure of $C_0(B)\cap \scr{F}$ in $\scr{F}$,   %$C_0(\mathbb{R})$ is dense in $\scr{F}$,  
we have
\begin{equation}\label{core}
	\begin{split}
		\lambda_{0}(B)&=\inf \{\mathscr{E}(f, f): f\in C_0(B)\cap \scr{F}, \pi(f^{2})=1 \}\\
		&=\inf \{\mathscr{E}(f, f): f \in \mathscr{F}, \pi(f^{2})=1 \text { and }\left.f\right|_{B^{c}}=0\}.
	\end{split}
\end{equation}
We call $\lambda_0(B)$  {\it the local Dirichlet eigenvalue} on $B$.

Let %$U^B$ the Green function of $Y$ on $B$, and  
$U_X^B$ be the Green operator of $X$ on $B$ and $G_X^B(\cdot,\cdot)$
be the  Green (density) function of $X$, i.e. for any $x,y\in \mathbb{R}$,  
$U_X^B(x,\d y)=G_X^B(x, y)\d y$. Denote by  $\tau_B^X$  the first exit time from $B$ and define the addictive functional  $A_t\coloneqq\int_{0}^t\sigma(X_s)^{-\alpha}\d s$.
Then by  the basic transform formula for time-change (see \cite[Lemma A.3.7]{CM12}), we have for any $K\in\scr{B}(\mathbb{R})$,
\begin{equation}\label{timechange}
	\begin{aligned}
		U^B(x, K) &=\mathbb{E}_{x}\left[\int_{0}^{\tau_B} \mathbf{1}_{K}\left(Y_{t}\right) \mathrm{d} t\right]=\mathbb{E}_{x}\left[\int_{0}^{\infty}\mathbf{1}_{\left\{Y_{t} \in K, t<\tau_{B}\right\}} \mathrm{d} t\right] \\
		&=\mathbb{E}_{x}\left[\int_{0}^{\infty} \mathbf{1}_{\left\{X_{t} \in K, t<\tau_{B}^{X}\right\}} \mathrm{d} A_{t}\right] \\
		&=U_X^B\left(\sigma^{-\alpha} \mathbf{1}_{K}\right)(x).
	\end{aligned}
\end{equation}
Hence for any $f\in \mathscr{B}(\mathbb{R})$ with  $U^B|f|<\infty$, 
$$
\begin{aligned}
	U^Bf(x)
	&=U_X^B\left(\sigma^{-\alpha} f\right)(x)=\int_K f(y)	G_X^B(x, y)\sigma(y)^{-\alpha}\d y.
\end{aligned}
$$

Since $G^B(\cdot,\cdot)$  can be represented by the Green function of $X$,  the estimate of $G^B(\cdot,\cdot)$ is obtained from $G_X^B(\cdot,\cdot)$. 

For the  one-dimensional symmetric $\alpha$-stable process $X$ with $\alpha\in(1,2)$,  its  Green function $G_X^B(\cdot, \cdot)$ for a open set $B$ can be expressed explicitly. For example:

(1) (\cite[Lemma 4]{bz06}) $B=\mathbb{R} \backslash\{0\}$ :
\begin{equation}\label{Green0}
	G_{X}^{\{0\}^c}(x, y)=\frac{\omega_\alpha}{2}\left(|y|^{\alpha-1}+|x|^{\alpha-1}-|y-x|^{\alpha-1}\right),
\end{equation}
where $\omega_\alpha$ is given by \eqref{omegaalpha}.

(2) (\cite[(11)]{KAE20}) $B=[-1,1]^{c}:$
\begin{equation}\label{Green1}
	G_X^{[-1,1]^c}(x,y)=c_{\alpha}\left(|x-y|^{\alpha-1}h\left(\frac{|xy-1|}{|x-y|}\right)-(\alpha-1)h(x)h(y)\right),
\end{equation}
where 	$c_{\alpha}=2^{1-\alpha} /(\Gamma(\alpha/2)^2)$, and
%$h$ is the harmonic function for $P_{t}^{[-1,1]^{c}}:$
\begin{equation}\label{h}
	h(x)=\int_{1}^{|x|}(z^{2}-1)^{\frac{\alpha}{2}-1} \mathrm{~d} z
\end{equation}
is the {\it harmonic function} for   $P_t^{[-1,1]^c}$, i.e. $P_t^{[-1,1]^c}h(x)=h(x)$ for any $x\notin [-1,1]$ (see \cite[Theorem 1.2]{KAE20}).

(3) (\cite[Page 388]{BB01}) $B=(0, \infty)$: 
\begin{equation}\label{Green_0infty}
	G_{X}^{(0, \infty)}(x, y)=\frac{1}{\Gamma(\alpha / 2)^{2}}|x-y|^{\alpha-1} J_{\alpha}\left(\frac{x \wedge y}{|x-y|}\right),
\end{equation}
where $J_{\alpha}(t)\coloneqq\int_{0}^{t}[s(s+1)]^{\frac{\alpha}{2}-1} \mathrm{~d} s.$

\section{Exponential ergodicity }\label{proof}

It is well known that the exponential ergodicity for a reversible Markov process is equivalent to the existence of the spectral gap $\lambda_{1}$, %In general, it is difficult to estimate $\lambda_1$ for L\'{e}vy-type jump processes directly, but
and we can turn this problem to the estimate of the local Dirichlet eigenvalue $\lambda_{0}(B)$ for some open set $B$.

By \cite[Theorem 1.4]{cmf00} or \cite[Theorem 1.1]{wfy99}, we have an upper bound for $\lambda_{1}$ by using $\lambda_{0}(B)$:%the (local) first Dirichlet eigenvalue:
\begin{equation}\label{upperbound}
	\lambda_{1}\leqslant\frac{\lambda_{0}(B)}{\pi(B^c)},\ \text{for any open set} \ B \ \text{with} \ \pi(B^c)>0,
\end{equation}
and by \cite[Proposition 3.2]{cmf00'}, $\lambda_{1}\geqslant\lambda_{0}\coloneqq\lambda_{0}(\{0\}^c)$. 
Therefore, to prove exponential ergodicity, our strategy is to estimate the local Dirichlet eigenvalues. 

To estimate the upper bound and lower bound of $\lambda_{0}(B)$, %we carefully choose a proper test function in the definition \eqref{vari diri}.
%For the lower bound, 
according to   the definition \eqref{vari diri} and \cite[Theorem 3.2]{st05}, we have the following variational formula for  $\lambda_{0}(B)$. %(more details will be discussed in the proof of Theorem \ref{Diri eigen}).
\begin{lemma}[Variational formula for the local Dirichlet eigenvalue]\label{variational}
	Assume that $B$ is a nonempty open subset of $\mathbb{R}$. %regular (see \cite[Page 68]{chung86} for the definition). 
	Then 
	$$
	\inf_{f\in C_b(\mathbb{R})}
	\sup_{x \in B} \frac{f(x)}{U^B f(x)}	\geqslant	\lambda_{0}(B)  \geqslant \sup_{f\in C_b(B)}
	\inf _{x \in B} \frac{f(x)}{U^B f(x)}, 
	$$
	where  $C_b(B)$ is the space of all bounded continuous functions on $B$.
\end{lemma}
\prf
First we consider the upper bound. Note that for $f\in C_b(\mathbb{R})$,  $\int_{\mathbb{R}}|f(x)|U^B|f|(x) \pi(\d x)<\infty$. By \cite[Theorem 1.3.9]{OY13}, we see that	
$U^Bf\in\mathscr{F}^B$, and 
\eqref{poisson weak-solution} holds.  
Thus by the definition \eqref{vari diri},
$$\lambda_{0}\leqslant \frac{\scr{E}(U^Bf,U^Bf)}{\pi((U^Bf)^2)}=\frac{\int f U^Bf\d \pi}{\pi((U^Bf)^2)}\leqslant \sup_{x \in B} \frac{f(x)}{U^B f(x)}.$$
So we get the upper bound by the arbitrariness of $f\in C_b(\mathbb{R}).$

For the lower bound, %since the semigroup $P_t$ of $Y$ is Feller and $P_t{1}\in C_b(\mathbb{R})$, we know that $P_t$ is a $C_b$-Feller semigroup (see \cite[Theorem 1.9]{bsw13}). By \cite[Proof of Theorem 1.7]{CW14}, the transition kernel $P_t (x, \d y)$ is absolutely continuous with respect to $\pi$, i.e. the transition density (heat kernel) $p_t(x,y)$ exists. Thus according to \cite[Theorem 1.14]{bsw13}, $Y$ is a strong Feller process. Therefore, due to \cite[Page 68]{chung86}, if $B$ is regular,  then the semigroup $P_t^{B}$ associated with the killed process $Y^B$ is a Feller process on $C_0(B)$. 
by \cite[Theorem 3.5.7(ii)]{OY13}, $Y^B$ is also a Hunt process, thus it is a right continuous Markov process. 
Hence by the proof of \cite[Lemma 2.2]{st05}, for any $f\in C_b(B)$,
\begin{equation}\label{resol}
	\widetilde{L}^B U_\beta^B f(x)=\beta U_\beta^B f(x)-f(x),
\end{equation}
where  $U_\beta^B f\coloneqq\int_{0}^\infty \e^{-\beta t}P_t^Bf(x)\d t$ and $\widetilde{L}^B$ is the weak generator for $Y^B$ (see \cite[Definition 2.1]{st05}). Therefore,  by \cite[Theorem 3.2]{st05} and \eqref{resol}, %we can estimate the local Dirichlet eigenvalues:
\begin{equation*}\label{vari}
	\lambda_{0}(B) \geqslant \inf _{x \in B}\left(-\frac{\widetilde{L}^B U_\beta^B f}{U_\beta^B f}\right)(x)=\inf _{x \in B} \frac{f(x)}{U_\beta^{B} f(x)}-\beta\geqslant \inf _{x \in B} \frac{f(x)}{U^{B} f(x)}-\beta.
\end{equation*}
Now 
we obtain the lower bound by letting $\beta \rightarrow 0$.
\deprf

\subsection{Estimates for the local Dirichlet eigenvalues}
In this section, our main aim is to  estimate the  bounds of the local Dirichlet eigenvalues  on  $\mathbb{R}\setminus\{0\}$ and $(0,\infty)$ by using \eqref{vari diri} and Lemma \ref{variational}.  

Before stating the main results, we recall the so-called \roman2-operator 
\begin{equation}\label{2-operator}
	\mathrm{\roman2} (f)(x)\coloneqq \frac{1}{f(x)}U_0f(x), \quad  f\in\{g: g\in C([0,\infty]), g(0)=0, g>0\} 
\end{equation}
constructed by M.F. Chen (cf. \cite[Section 6.2]{cmf05}) for  diffusion operator $$\mathcal{A}=a(x)\frac{\d^2}{\d x^2}+b(x)\frac{\d}{\d x}$$ on half line,
where  $a>0$ and $b$ are continuous on $[0,\infty)$, $$U_0f(x)\coloneqq\int_0^x\e^{-C(y)}\left(\int_y^\infty f(z)\nu(\d z)\right)\d y,$$  $C(x)=\int_1^x(b(t)/a(t))\d t$ and $\nu(\d x)=(\e^{C(x)}/{a(x)}) \d x$ is the reversible measure. 
By choosing $f=\sqrt{\varphi}\coloneqq\sqrt{\int_0^x\e^{-C(t)}\d t}$ in \eqref{2-operator}, \cite[Theorem 6.1]{cmf05} gives the explicit estimate for the Dirichlet eigenvalue
$$\lambda_{0}\geqslant  (\inf_x \mathrm{\roman2} (\sqrt{\varphi})(x))^{-1}\geqslant (4\eta)^{-1},$$
where $\eta\coloneqq\sup_{x>0}\varphi(x)\pi([x,\infty))$. 

We interpret the \roman2-operator from a new viewpoint. It is well known that $U_0f(x)$ is the classical solution  of ordinary differential equation $-\mathcal{A}(U_0f)=f$ with Dirichlet boundary condition  $U_0f(0)=0$, and  $\varphi(x)=\int_0^x\e^{-C(t)}\d t$ is the harmonic function of $\mathcal{A}$ with $\varphi(0)=0$. By choosing the square root of harmonic function in the \roman2-operator, we obtain the lower bound for $\lambda_0$ of diffusion operator. 

%Furthermore,  the  optimal test function is the square root of harmonic function $\varphi$. Therefore,
%choosing $f(x)=\sqrt{\varphi}$, we  obtain that
%$$\lambda_0\geqslant(4\delta)^{-1}.$$

Now  we   construct the \roman2-operator for time-changed symmetric stable process $Y$.  
Let 
$$U^{(0)}\coloneqq U^{\{0\}^c}$$
be the Green operator of $Y$ killed upon $\{0\}^c$,  by \eqref{poisson weak-solution}, $v\coloneqq U^{(0)}f$ is a weak solution of $-Lv=f$ on $L^2(\pi)$ with $v(0)=0$. 

Naturally, we define the \roman2-operator by
$$\mathrm{\roman2}(f)(x)\coloneqq\frac{1}{f(x)}U^{(0)}f(x),\quad \text{for}\quad f\in\mathscr{G}\coloneqq\{g: g\in C(\mathbb{R}), U^{(0)}|g|<\infty, g(0)=0, g>0\}.$$
%By a basic estimate, it has the following representation which is similar with \eqref{2-operator}:
\begin{lemma} For any $x\in\mathbb{R}$ and $f\in \mathscr{G}$, 
	$$\mathrm{\roman2}(f)(x)\leqslant\frac{(\alpha-1)\omega_\alpha}{f(x)}\int_{0}^{|x|} z^{\alpha-2}\left(\int_{\mathbb{R}\setminus(-z,z)} f(y)\sigma(y)^{-\alpha}\d y\right)\d z.$$
\end{lemma}
\prf
By the property of time change \eqref{timechange}, for any $f$ with $U^{(0)}|f|<\infty$,  
$$U^{(0)}f(x)=\int_{\mathbb{R}}G_X^{\{0\}^c}(x,y)  f(y)\sigma(y)^{-\alpha}\d y.$$
By an elementary inequality (see \cite[Lemma 4.2]{mwt20}): 
\begin{equation}\label{basic ineq}
	|y|^{\alpha-1}+|x|^{\alpha-1}-|y-x|^{\alpha-1}\leqslant2(|x|\wedge|y|)^{\alpha-1}, \  1<\alpha<2,
\end{equation}
we obtain that for any $x,y\neq 0$,
\begin{equation*}\label{esti for G0}
	G_{X}^{\{0\}^c}(x, y)=\frac{\omega_\alpha}{2}\left(|y|^{\alpha-1}+|x|^{\alpha-1}-|y-x|^{\alpha-1}\right)\leqslant \omega_\alpha(|y|\wedge |x|)^{\alpha-1}.
\end{equation*}
Then 
for any $f$ with $	U^{(0)} |f|<\infty$,
\begin{equation*}
	\begin{split}
		U^{(0)}f(x)
		&\leqslant \omega_\alpha\int_{\mathbb{R}}(|y|\wedge |x|)^{\alpha-1}f(y)\sigma(y)^{-\alpha}\d y\\
		&=\omega_\alpha\left(\int_{\mathbb{R}\setminus(-|x|,|x|)} |x|^{\alpha-1}f(y)\sigma(y)^{-\alpha}\d y+\int_{-|x|}^{|x|} |y|^{\alpha-1}f(y)\sigma(y)^{-\alpha}\d y\right)\\
		&=(\alpha-1)\omega_\alpha\int_{0}^{|x|} z^{\alpha-2}\left(\int_{\mathbb{R}\setminus(-z,z)} f(y)\sigma(y)^{-\alpha}\d y\right)\d z.\\
	\end{split}
\end{equation*}
\deprf

Let $h_0(x)=(\omega_\alpha/2)|x|^{\alpha-1}$ be the harmonic function for $P_t^{\{0\}^c}$ on $\mathbb{R}\setminus\{0\}$ (see  \cite[Example 1.1]{ya10}). By choosing $f(x)=\sqrt{h_0(x)}, $ we have the following upper estimate:
\begin{lemma}\label{delta}
	If $$
	\delta\coloneqq\sup _{x} |x|^{\alpha-1} \int_{\mathbb{R} \backslash(-|x|, |x|)} \sigma(z)^{-\alpha} \mathrm{d} z<\infty,
	$$ then for any $x\in\mathbb{R}$,
	\begin{equation*}
		\mathrm{\roman2}(\sqrt{h_0})(x)\leqslant4\omega_{\alpha}\delta.
	\end{equation*}	
	%Specially, for symmetric case $\rho=1/2$, $$-\frac{4A_\alpha\delta}{\pi}=-\frac{4\delta}{\Gamma(\alpha)\cos(\pi\alpha/2)}.$$
\end{lemma}
\prf Note for any $y>0$,
\begin{equation*}\label{int by part}
	\int_{\mathbb{R}\setminus(-y,y)} |z|
	^{(\alpha-1)/2}\sigma(z)^{-\alpha}\d z
	=	\int_{y}^\infty z
	^{(\alpha-1)/2}\sigma(z)^{-\alpha}\d z+\int_{-\infty}^{-y} (-z)
	^{(\alpha-1)/2}\sigma(z)^{-\alpha}\d z.
\end{equation*}
By using integration by parts, for any $y>0$, we have
\begin{equation*}
	\begin{split}
		\int_{y}^\infty z
		^{(\alpha-1)/2}\sigma(z)^{-\alpha}\d z&=-\int_{y}^\infty z
		^{(\alpha-1)/2}\d \mu((z,\infty))\\
		&\leqslant 	y
		^{(\alpha-1)/2}\mu((y,\infty))+\frac{\alpha-1}{2}\int_{y}^{\infty}z
		^{(\alpha-3)/2}\mu((z,\infty))\d z,\\
	\end{split}
\end{equation*}
while
\begin{equation*}
	\begin{split}
		\int_{-\infty}^{-y} (-z)
		^{(\alpha-1)/2}\mu(\d z)&=\int_{y}^\infty z
		^{(\alpha-1)/2}\sigma(-z)^{-\alpha}\d z=-\int_{y}^\infty z
		^{(\alpha-1)/2}\d \mu((-\infty,-z))\\
		&\leqslant 	y
		^{(\alpha-1)/2}\mu((-\infty,-y))+\frac{\alpha-1}{2}\int_{y}^{\infty}z
		^{(\alpha-3)/2}\mu((-\infty,-z))\d z.\\
	\end{split}
\end{equation*}
Therefore,
\begin{equation*}\label{int by part}
	\int_{\mathbb{R}\setminus(-y,y)} |z|
	^{(\alpha-1)/2}\mu(\d z)
	\leqslant 	y
	^{(\alpha-1)/2}\mu((-y,y)^c)+\frac{\alpha-1}{2}\int_{y}^{\infty}z
	^{(\alpha-3)/2}\mu((-z,z)^c)\d z.
\end{equation*}
Note the definition of  $\delta$ gives  
\begin{equation*}
	\begin{split}
		\int_{\mathbb{R}\setminus(-y,y)} |z|
		^{(\alpha-1)/2}\mu(\d z)
		\leqslant\frac{\delta}{y^{(\alpha-1)/2}}+\frac{\delta(\alpha-1)}{2}\int_{y}^{\infty}z^{-(\alpha+1)/2}\d z%+\varphi(n)\mu(({n},{\infty}))
		=\frac{2\delta}{y^{(\alpha-1)/2}}.\\
	\end{split}
\end{equation*}
Finally we have  \begin{equation*}
	\begin{split}
		\mathrm{\roman2}(\sqrt{h})(x)&\leqslant\frac{(\alpha-1)\omega_\alpha}{|x|^{(\alpha-1)/2}}\int_{0}^{|x|} z^{\alpha-2}\frac{2\delta}{z^{(\alpha-1)/2}}\d z={4\omega_\alpha\delta}.
	\end{split}
\end{equation*}	
\deprf

Next, we  consider the local  Dirichlet eigenvalue on half-line. 
According to \cite[Lemma 5]{BB01},
$$
\frac{(\alpha-1) \Gamma(\alpha / 2)^{2} G_{X}^{(0, \infty)}(x, y)}{(x \wedge y)^{\alpha-1}} \leqslant 1.
$$
For any $f$ with $	U^{(0, \infty)} |f|<\infty$, it follows from  \eqref{timechange} that
%$$U^B f(x)=U^B^{X}\left(f \sigma^{-\alpha}\right)(x), x \in \mathbb{R}.$$
$$
\begin{aligned}
	U^{(0, \infty)} f(x) &=U_{X}^{(0, \infty)}\left(f \sigma^{-\alpha}\right)(x)=\int_{0}^{\infty} G_{X}^{(0, \infty)}(x, y) f(y) \sigma(y)^{-\alpha}\d y \\
	& \leqslant \frac{1}{(\alpha-1) \Gamma(\alpha / 2)^{2}} \int_{0}^{\infty}(x \wedge y)^{\alpha-1} f(y) \sigma(y)^{-\alpha}\d y \\
	&=\frac{1}{(\alpha-1) \Gamma(\alpha / 2)^{2}}\left(\int_{0}^{x} y^{\alpha-1} f(y) \sigma(y)^{-\alpha}\d y+\int_{x}^{\infty} x^{\alpha-1} f(y) \sigma(y)^{-\alpha}\d y\right) \\
	&=\frac{1}{\Gamma(\alpha / 2)^{2}} \int_{0}^{x} z^{\alpha-2}\left(\int_{z}^{\infty} f(y) \sigma(y)^{-\alpha}\d y\right) \mathrm{d} z.
\end{aligned}
$$
By  defining $$
\begin{aligned}
	\mathrm{\roman2}^{+}(f)(x) &\coloneqq\frac{1}{f(x)} U^{(0, \infty)} f(x)=\frac{1}{\Gamma(\alpha / 2)^{2}f(x)} \int_{0}^{x} z^{\alpha-2}\left(\int_{z}^{\infty} f(y) \sigma(y)^{-\alpha}\d y\right) \mathrm{d} z,\\
	%I I^{-}(f)(x) &=\frac{1}{f(x)} G_{(-\infty, 0)} f(x)=\frac{c_{1}(\alpha)}{f(x)} \int_{0}^{x} z^{\alpha-2}\left(\int_{-\infty}^{-z} f(y) \sigma(y)^{-\alpha}\d y\right) \mathrm{d} z
\end{aligned}
$$
%Let $\varphi(x)=x^{(\alpha-1) / 2}$. 
a similar argument to Lemma \ref{delta} gives the following estimate.%to the proof of Lemma \ref{II+}, we have the following estimate, and we omit the details here.
\begin{lemma}\label{delta+}If
	$$\delta_{+}\coloneqq\sup _{x>0} x^{\alpha-1} \int_{x}^{\infty} \sigma(z)^{-\alpha} \mathrm{d} z<\infty,$$
	then
	$$\mathrm{I I}^{+}(\varphi)(x) \leqslant \frac{4}{\Gamma(\alpha / 2)^{2}(\alpha-1)} \delta_{+},$$
	where $\varphi(x)=x^{(\alpha-1) / 2}$.

\end{lemma}

To obtain the lower bounds of the local  Dirichlet eigenvalues, we use the following approximation, which is modified from \cite[Lemma 3.4]{cmf00''} for jump process.
\begin{lemma}\label{approach}
	Let $A$ be an  open subset of $\mathbb{R}$, and $\{A_m\}_{m=1}^{\infty}$ be a sequence of bounded open subsets such that $A_m \uparrow A$. Then we have %there exists a subsequence $A_{m_n} \uparrow A,$ such that 
	$$\lambda_{0}(A)=\lim _{m\rightarrow \infty} \lambda_{0}(A_{m}).
	$$
\end{lemma}
\prf 
%Assume that $A$ is unbounded, otherwise, we let $A_m\equiv A$ for any $m\geqslant 1$. 
Since $A_{m}\subset A$, by the definition of the local Dirichlet eigenvalue \eqref{vari}, we have $\lambda_{0}(A)\leqslant\lambda_{0}(A_{m})$. So we only need to prove
$\lambda_{0}(A)\geqslant\lim _{m \rightarrow \infty} \lambda_{0}(A_{m})$.

%Since $C_0^\infty(\mathbb{R})$ which is the space of all smooth functions with compact support on $\mathbb{R}$  is a core of $\scr{F}$, 
By  \eqref{core}, for every $m\geqslant 1,$ there is $f_{m} \in C_0(A)\cap \scr{F}$,  $\pi\left(f_{m}^{2}\right)=1$,   such that 
\begin{equation}\label{A-K_n}
	\lambda_{0}(A) \geqslant \mathscr{E}\left(f_{m},f_m\right)-\frac{1}{m} \geqslant \lambda_{0}(K_m)-\frac{1}{m},
\end{equation}
where $K_m\coloneqq\supp{f_m}\subset A$ is the compact support of $f_m$. %Since $K_n$ is bounded, $K_n\subset A$ follows from the unboundedness of $A$.  

Since $A_m \uparrow A,$  for each $m$, there exists 
$k_m$ which satisfies that $k_m\uparrow\infty$ as $m\rightarrow\infty$, such that %for any $m\geqslant k_m$, 
$K_m\subset A_{k_m}$, so then $\lambda_{0}(A_{k_m})\leqslant \lambda_{0}(K_m)$. By combining it with \eqref{A-K_n},  we have
$$\lambda_{0}(A) \geqslant  \lambda_{0}(A_{k_m})-\frac{1}{m}.$$
Due to the monotonicity of $\{A_m\}$,   the required assertion follows by   letting  $m\rightarrow\infty$.

\deprf

%We are now in the positive to prove  Theorem \ref{Diri eigen}.

Now we establish the  estimates for the local Dirichlet eigenvalues.

\noindent\textbf{Proof of Theorem \ref{Diri eigen}}.  
Let $A=(0,\infty)$ (or $\mathbb{R}\setminus\{0\}$) and  $A_n=(-n,n)\cap A$, for $n\geqslant1$. Denote by  $Y_{t}^{n}$ the killed process on $A_n$. %and $L_{n}$ be the weak generator of  $Y^{n}$.  

Since %$A_n$ is regular, and 
the continuous function  is bounded on $A_n$, by Lemma \ref{variational}, we have for any $f\in C(A_n)$, 
$$\lambda_{0}(A_n) \geqslant %\inf _{x \in A_n}\left(-\frac{\widetilde{L}_{n} G_{\beta}^{A_n} f(x)}{G_{\beta}^{A_n} f(x)}\right) = 
\inf _{x \in A_n} \frac{f(x)}{U^{A_n} f(x)}.$$
It is clear that $U^{A_n} f\leqslant U^{A} f$ by noting  $\tau_{A_n}\leqslant\tau_{A}$.
By Lemma \ref{approach}, % there exists a subsequence $\{m_n\} \uparrow \infty,$ such that 
$\lambda_{0}(A)=\lim _{n\rightarrow \infty} \lambda_{0}(A_{n})$, so 
$$
\begin{aligned}
	\lambda_{0}(A) &=\lim _{n \rightarrow \infty} \lambda_{0}(A_n) \geqslant \lim _{n \rightarrow \infty} \inf _{x \in A_n} \frac{f(x)}{U^{A} f(x)}. \\
\end{aligned}
$$
In the case of $A=(0,\infty)$ and $A_n=(0,n)$,   by letting $f(x)=x^{(\alpha-1)/2}$, we have
$$
\begin{aligned}
	\lambda_{0}((0, \infty)) & \geqslant \lim _{n\rightarrow \infty} \inf _{x \in(0,n)} \frac{f(x)}{U^{(0,\infty)} f(x)}  \geqslant \inf \left(\mathrm{I I}^{+}(f)(x)\right)^{-1} \\
	& \geqslant\left[\frac{4}{\Gamma(\alpha / 2)^{2}(\alpha-1)} \delta\right]^{-1}=\frac{(\alpha-1)\Gamma(\alpha / 2)^{2}}{4\delta};
\end{aligned}
$$
in the case of $A=\mathbb{R}\setminus\{0\}$,  and $A_n=(-n,n)\setminus\{0\}$,   by letting $f(x)=|x|^{(\alpha-1)/2}$, we have
$$
\begin{aligned}
	\lambda_{0} & \geqslant \lim _{n \rightarrow \infty} \inf _{x \in(0, n)} \frac{f(x)}{U^{(0)}f(x)}  \geqslant \inf \left(\mathrm{\roman2}(f)(x)\right)^{-1}  \geqslant\frac{1}{4\omega_\alpha\delta}.
\end{aligned}
$$

Next we estimate the upper bound. For this, we assume that $\lambda_{0}>0$, otherwise, the conclusion is trivial.  

%Let  $\{E_{\lambda}^0\}$ be the spectral measure of $L^{\mathbb{R}\setminus\{0\}}$. Note that $\int_{(0,\lambda_{0})}\d E_{\lambda}^0=0$, then we have for any $g\in L^2(\mathbb{R}\setminus\{0\})$,
Since
\begin{equation*}\label{spec decom}
	\|P_t^{\{0\}^c}g\|_{L^2(\pi)}%=\int_{\mathbb{R}}\left(\int_{0}^\infty\e^{-\lambda t}\d E_{\lambda}g\right)^2\d \pi
	\leqslant \e^{-\lambda_0 t}\|g\|_{L^2(\pi)},
\end{equation*}
we have that
$$\|U^{(0)}g\|_{L^2(\pi)}\leqslant\lambda_{0}^{-1}\|g\|_{L^2(\pi)},$$
which implies that $(\scr{E}^{\{0\}^c}, \scr{F}^{\{0\}^c})$ is transient, and $$\int gU^{(0)}g\d\pi\leqslant \|g\|_{L^2(\pi)}\|U^{(0)}g\|_{L^2(\pi)}\leqslant\frac{1}{\lambda_0}\|g\|_{L^2(\pi)}^2<\infty.$$
By \cite[Theorem 1.3.9]{OY13}, we have $U^{(0)}g\in\scr{F}$, and \eqref{poisson weak-solution} holds. Hence
$$\scr{E}(U^{(0)}g,U^{(0)}g)=\int gU^{(0)}g\d \pi.$$
According to the definition of $\lambda_{0}$,
\begin{equation}\label{upbd}
	\lambda_{0}\leqslant\frac{\scr{E}(U^{(0)}g,U^{(0)}g)}{\|U^{(0)}g\|_{L^2(\pi)}^2}=\frac{\int gU^{(0)}g\d \pi}{\|U^{(0)}g\|_{L^2(\pi)}^2}\leqslant \sup_{x>0}\frac{g(x)}{U^{(0)}g(x)}+\sup_{x<0}\frac{g(x)}{U^{(0)}g(x)}.
\end{equation}
Note that for $xy>0$, $$|y|^{\alpha-1}+|x|^{\alpha-1}-|y-x|^{\alpha-1}\geqslant(|x|\wedge|y|)^{\alpha-1}, \  1<\alpha<2,$$
so for $x>0$,
\begin{equation}\label{positive}
	{U^{(0)}g(x)}\geqslant\frac{\omega_\alpha}{2}\int_{0}^{\infty}(x\wedge y)^{\alpha-1}g(y)\sigma(y)^{-\alpha}\d y, \quad%\geqslant\frac{\omega_\alpha}{2}x^{\alpha-1}\int_{x}^\infty g(y)\sigma(y)^{-\alpha}\d y.
\end{equation}
while for $x<0$,
\begin{equation}\label{negative}
	{U^{(0)}g(x)}\geqslant\frac{\omega_\alpha}{2}\int_{-\infty}^0(-(x\vee y))^{\alpha-1}g(y)\sigma(y)^{-\alpha}\d y.
\end{equation}
Fix $x_0>0$ and choose 
$$g(x)=\left\{\begin{array}{ll}
	(x\wedge x_0)^{\alpha-1}, & x \geqslant0, \\
	|x\vee(-x_0)|^{\alpha-1}, & x \leqslant0.
\end{array}\right.$$
Then $g\in L^2(\pi)$ and $g(0)=0$. For $x>0$, by \eqref{positive}, 
\begin{equation}\label{G+}
	\begin{split}
		\frac{U^{(0)}g(x)}{g(x)}&\geqslant\frac{\omega_\alpha}{2}\frac{1}{(x\wedge x_0)^{\alpha-1}}\int_{x\wedge x_0}^\infty(x\wedge y)^{\alpha-1}(y\wedge x_0)^{\alpha-1}\sigma(y)^{-\alpha}\d y\\
		&\geqslant \frac{\omega_\alpha}{2}\int_{x\wedge x_0}^\infty(y\wedge x_0)^{\alpha-1}\sigma(y)^{-\alpha}\d y\geqslant \frac{\omega_\alpha}{2}\int_{ x_0}^\infty(y\wedge x_0)^{\alpha-1}\sigma(y)^{-\alpha}\d y\\
		&\geqslant\frac{\omega_\alpha}{2}x_0^{\alpha-1}\int_{ x_0}^\infty\sigma(y)^{-\alpha}\d y.
	\end{split}
\end{equation}
For $x<0$, by \eqref{negative},
\begin{equation}\label{G-}
	\begin{split}
		\frac{U^{(0)}g(x)}{g(x)}&\geqslant\frac{\omega_\alpha}{2}\frac{1}{|x\vee (-x_0)|^{\alpha-1}}\int_{-\infty}^{x\vee (-x_0)}(-(x\vee y))^{\alpha-1}|y\vee (-x_0)|^{\alpha-1}\sigma(y)^{-\alpha}\d y\\
		&\geqslant \frac{\omega_\alpha}{2}\int_{-\infty}^{x\vee (-x_0)}|y\vee (-x_0)|^{\alpha-1}\sigma(y)^{-\alpha}\d y\geqslant \frac{\omega_\alpha}{2}\int_{-\infty}^{-x_0}|y\vee (-x_0)|^{\alpha-1}\sigma(y)^{-\alpha}\d y\\
		&\geqslant\frac{\omega_\alpha}{2}x_0^{\alpha-1}\int_{-\infty}^{-x_0}\sigma(y)^{-\alpha}\d y.
	\end{split}
\end{equation}
Now combining  \eqref{upbd}, \eqref{G+}, and \eqref{G-}, we obtain the desired result.

\deprf

\subsection{Proof of Theorem \ref{exp erg}}
By combining Theorem \ref{Diri eigen} with  \eqref{vari diri}, we prove the criterion for exponential ergodicity, and bounds for spectral gap $\lambda_{1}$.

\noindent\textbf{Proof of Theorem \ref{exp erg}}. 
The sufficiency and the estimate for lower bound of $\lambda_{1}$ follow from $\lambda_{1}\geqslant\lambda_{0}$ and Theorem \ref{Diri eigen}, so it remains to show the necessity and the upper bound. 

First, without loss of generality, assume that 
\begin{equation*}\label{delta+}
	\delta_{+}\coloneqq\sup _{x>0} x^{\alpha-1} \int_{x}^{\infty} \sigma(y)^{-\alpha} \mathrm{d} y=\infty.
\end{equation*}
Let $h$ be  the harmonic function for $P_t^{[-1,1]^c}$ given by \eqref{h}. Note that for any $x \in \mathbb{R},$ 
$$x^{\alpha-1} \int_{x}^{\infty} \sigma(y)^{-\alpha} \mathrm{d} y<\infty,\ \text{and}\
\lim _{x \rightarrow \infty} \frac{h(x)}{x^{\alpha-1}}=\frac{1}{\alpha-1}.
$$
So $\delta_{+}=\infty$ means that
\begin{equation}\label{delta+}
	\delta_{+}(x)\coloneqq h(x) \int_{x}^{\infty} \sigma(y)^{-\alpha} \mathrm{d} y \rightarrow \infty, \text { as } x \rightarrow \infty.
\end{equation}
By \cite[ Lemma 3.3]{KAE20}, for any $x \notin[-1,1]$,
\begin{equation}\label{limit}
	\lim _{y \rightarrow \infty} G_{X}^{[-1,1]^{c}}(x, y)=K_{\alpha} h(x),
\end{equation}
where  $K_{\alpha}$  are given by  $$K_{\alpha}=\frac{2 c_{\alpha}\left(1-\frac{\alpha}{2}\right) \Gamma\left(\frac{\alpha}{2}\right)}{\Gamma\left(1-\frac{\alpha}{2}\right)} \int_{1}^{\infty} \frac{h^{\prime}(v)}{1+v} \mathrm{~d} v<\infty.$$
Thus there exists some constant $N_1>1,$ such that for any $y>N_1$,
$$
G_{X}^{[-1,1]^{c}}(x, y) \geqslant \frac{1}{2} K_{\alpha} h(x), \quad x\notin [-1,1].
$$
For any fixed $x_0>N_1$, let $h^{x_0}(x)\coloneqq h(x\wedge x_0)$, and $u^{x_0}(x)\coloneqq\left(U^{[-1,1]^{c}} h^{x_0}\right)\left(x \right)$. Note that
$$
\begin{aligned}
	u^{x_0}(x) &=U_{X}^{[-1,1]^{c}}\left(h^{x_0} \sigma^{-\alpha}\right)\left(x \right) \geqslant \int_{x_0}^{\infty} G_{X}^{[-1,1]^{c}}(x, y) h^{x_0}(y) \mu(\mathrm{d} y) \\
	& \geqslant \frac{K_\alpha}{2}h\left(x \right) h\left(x_{0}\right) \int_{x_{0}}^{\infty} \mu(\mathrm{d} y)=\frac{K_\alpha}{2}h\left(x \right) \delta_{+}\left(x_{0}\right)\geqslant\frac{K_\alpha}{2}h^{x_0}\left(x \right) \delta_{+}\left(x_{0}\right).
\end{aligned}
$$
Hence
%$$u^{x_{0}}\left(x \wedge x_{0}\right) / h\left(x \wedge x_{0}\right) \geqslant \delta_{+}\left(x_{0}\right).$$ 
$$\langle u^{x_{0}}, h^{ x_{0}} \rangle_{\pi}=\left\langle \left(u^{x_{0}}\right)^{2}, \frac{h^{ x_{0}}} {u^{x_{0}}}\right\rangle_\pi\leqslant \|u^{x_{0}}\|_{L^2(\pi)}^2\left(\frac{K_\alpha}{2}\delta_{+}\left(x_{0}\right)\right)^{-1}<\infty,$$ 
where $\lan\cdot,\cdot\ran_\pi$ is the inner product on $L^2(\pi)$.  By \cite[Theorem 1.3.9]{OY13}, we have $u^{x_{0}}\in \scr{F}$ and  $\mathscr{E}\left(u^{x_{0}}, u^{x_{0}}\right)=\langle h^{x_{0}},u^{x_{0}}\rangle_{\pi}.$ 
According to  the definition of the local	Dirichlet eigenvalue,
\begin{equation}\label{K-upper}
	\begin{aligned}
		\lambda_{0}\left([-1,1]^{c}\right) & \leqslant \frac{\mathscr{E}\left(u^{x_{0}}, u^{x_{0}}\right)}{\|u^{x_{0}}\|_{L^2(\pi)}^2}=\frac{\left\langle u^{x_{0}}, h^{x_{0}}\right\rangle_{\pi}}{\|u^{x_{0}}\|_{L^2(\pi)}^2 }
		%	& \leqslant \frac{\left\langle h\left(\cdot \wedge x_{0}\right) / u^{x_{0}},\left(u^{x_{0}}\right)^{2}\right\rangle_\pi}{\pi\left(\left(u^{x_{0}}\right)^{2}\right)} 
		\leqslant\left(\frac{K_\alpha}{2}\delta_{+}\left(x_{0}\right)\right)^{-1}.
	\end{aligned}
\end{equation}

Therefore,  by letting $x_{0} \rightarrow \infty,$  from \eqref{delta+}, we obtain that
$\lambda_{0}\left([-1,1]^{c}\right)=0$ which implies that $Y$ is non-exponentially ergodic.
\deprf

\section{Strong ergodicity}

To prove the strong ergodicity of $Y$, we need to estimate the uniform bounds for the first moment of hitting time. 
%For the necessity,  by a similar argument to \cite[Lemma 2.1]{myh02}, we know that strong ergodicity implies that for any closed set $B\subset\mathbb{R}$ with $\pi(B)>0$, $\sup_x\mathbb{E}_x\tau_B<\infty$. For the sufficiency, we only need to prove $M_0\coloneqq\sup_x\mathbb{E}_x\tau_{\{0\}^c}<\infty$.  

\noindent {\bf Proof of Theorem \ref{ergodc}}.

By a similar argument to \cite[Lemma 2.1]{myh02}, we know that strong ergodicity implies that for any closed set $B\subset\mathbb{R}$ with $\pi(B)>0$, $\sup_x\mathbb{E}_x\tau_B<\infty$. So for the necessity, we only need to prove  $	\sup _{x} \mathbb{E}_{x} \tau_{[-1,1]^c}=\infty$ under the assumption  $I=\infty$.
According to \eqref{limit} and the symmetry of $G_{X}^{[-1,1]^{c}}(x, y)$, we have
\begin{equation}\label{limgreen}
	\lim _{x \rightarrow \infty} G_{X}^{[-1,1]^{c}}(x, y)=K_{\alpha} h(y), \\
\end{equation}
where $K_\alpha$ is defined by $$K_{\alpha}=\frac{2 c_{\alpha}\left(1-\frac{\alpha}{2}\right) \Gamma\left(\frac{\alpha}{2}\right)}{\Gamma\left(1-\frac{\alpha}{2}\right)} \int_{1}^{\infty} \frac{h^{\prime}(v)}{1+v} \mathrm{~d} v<\infty.$$
Therefore, by Fatou's lemma and \eqref{limgreen}, 
$$
\begin{aligned}
	\sup _{x} \mathbb{E}_{x} \tau_{[-1,1]^c} &=\sup _{x} \int_{\mathbb{R} \backslash[-1,1]} G_{X}^{[-1,1]^{c}}(x, y) \sigma(y)^{-\alpha} \mathrm{d} y \\
	& \geqslant \int_{1}^{\infty} \liminf_{x \rightarrow \infty} G_{X}^{[-1,1]^{c}}(x, y) \sigma(y)^{-\alpha} \mathrm{d} y \\
	&=\int_{1}^{\infty} K_{\alpha} h(y) \sigma(y)^{-\alpha} \mathrm{d} y \\
	& \geqslant \frac{K_{\alpha}}{(\alpha-1)} \int_{1}^{\infty}\left(y^{\alpha-1}-1\right) \sigma(y)^{-\alpha} \mathrm{d} y=\infty.
\end{aligned}
$$

The sufficiency and the lower bound for convergence rate 
are proved by \cite[Theorem 4.3]{mwt20}. We survey the proof below for the readers' convenience. By \eqref{timechange},  we have
$$M_0\coloneqq\sup_x\mathbb{E}_x\tau_{\{0\}^c}=\sup_x\int_{\mathbb{R}}U^{(0)}(x,\d y)=\sup_x\int_{\mathbb{R}}G_X^{\{0\}^c}(x,y) \sigma(y)^{-\alpha}\d y.$$
According to \eqref{basic ineq},  we have
\begin{equation}\label{M0}
	M_0\leqslant \omega_{\alpha}\int_{\mathbb{R}}|y|^{\alpha-1}  \sigma(y)^{-\alpha}\d y=\omega_{\alpha}I.
\end{equation}
Note that 
$\lambda_1\geqslant\lambda_0,$ and  $\lambda_{0}^{-1}\leqslant M_0$ (see \cite[Lemma 3.2]{GT12}). Therefore,
by  \cite[Theorem 1.2 (R2)]{mwt20}, $\kappa\geqslant\min\{\lambda_{1},M_0^{-1}\}=M_0^{-1}\geqslant(\omega_{\alpha}I)^{-1}$, so $Y$ is strongly ergodic.

\section*{Acknowledgements}
The author thanks Prof. Yong-Hua Mao for valuable conversations of this paper.  This work is supported in part by the  National Key Research and Development Program of China (2020YFA0712900), the National Natural Science Foundation of China (Grant No.11771047) and the project from the Ministry of Education in China.

\bibliographystyle{plain}
\bibliography{tcstable}

\end{document}